\documentclass[12pt]{article}
\usepackage[latin1]{inputenc}
\usepackage[T1]{fontenc}
\usepackage[english]{babel}
\usepackage{multirow}
\usepackage[colorlinks=true,linkcolor=blue,urlcolor=red,citecolor=red]{hyperref}
\usepackage{theorem}
\usepackage{amsmath,amssymb}
\usepackage{mathrsfs}
\usepackage{dsfont}
\usepackage{makeidx}
\usepackage{bbold} 
\usepackage{xcolor}
\usepackage{fancybox}
\usepackage{fancyhdr}
\usepackage{textcomp}
\usepackage{varioref}  
\usepackage{chapterbib}

\usepackage{graphicx,subfigure}
\usepackage{color}
\usepackage{amstext}

\newcommand{\R}{\mathbb{R}}

\newcommand{\E}{\mathrm{E}}
 
\newcommand{\partiel}{\mathrm{d}}

\newcommand{\law}{\mathrm{Law}}

\newcommand{\var}{\mathrm{Var}}

 \theoremstyle{changebreak}
 \newtheorem{Theorem}{Theorem}[section]

 \theoremstyle{changebreak}
 \newtheorem{Definition}[Theorem]{Definition}
 \theoremstyle{changebreak}
 
 \theoremstyle{changebreak}
 \newtheorem{Lemma}[Theorem]{Lemma}

 \theoremstyle{changebreak}
 \newtheorem{proposition}[Theorem]{Proposition}

 \theoremstyle{changebreak}

 \theoremstyle{changebreak}
 \newtheorem{Corollary}[Theorem]{Corollary}

 \theoremstyle{changebreak}
 \newtheorem{Remark}[Theorem]{Remark} 

\fancypagestyle{main}{

\fancyhead[L]{\textsc{{\color{purple}\rule[2mm]{4mm}{1mm}}\hspace{0.1cm}}{\color{blue}\rule[2mm]{1mm}{1mm}}}

\fancyfoot[L]{\textsc{{\color{purple}\rule[2mm]{4mm}{1mm}}\hspace{0.1cm}}{\color{blue}\rule[2mm]{1mm}{1mm}}}
}

\author{Mack Dowell {\sc Komba Moudoumou}\footnote{\textbf{mackdowk@gmail.com}}\; , Fulgence {\sc Eyi Obiang}\footnote{eyiobiang@yahoo.fr}\; \\ and Octave {\sc Moutsinga}\footnote{octavemoutsing@gmail.com}}
\date{ {\small  Universit\'e des Sciences et Techniques de Masuku\\
		Facult\'e des Sciences - Dpt Math\'ematiques et Informatique\\
		Unit\'e de Recherches en Math\'ematiques et Informatique\\
		BP 943 Franceville, Gabon.}}
\title{\textbf{Generalized variational principle for an accelerated pressureless gas model}}
\begin{document}
\selectlanguage{english}
\makeatletter
\maketitle
\thispagestyle{main}
\begin{abstract}

We examine a discrete model of sticky particles initially subjected to acceleration. We propose a novel generalized variational principle for characterizing clusters (i.e., particle agglomerations) under decreasing acceleration function. 
This accelerated particle model is usually employed to solve a system of forced  (inhomogeneous) pressureless gas equations. The momentum conservation law requires evaluating the force acting on the gas particles. Here, we present a new solution to this system of equations, based on the distribution of mass over particle velocities and on the acceleration function.


\end{abstract}


\section{Introduction}

We are interested in the system of forced (inhomogeneous) pressureless gas equations on $\mathbb{R} \times \mathbb{R}_+$:
\begin{align}\label{equation gas disc}
\left\{
\begin{array}{ll}
\partial_{t}(\mu) + \partial_{x}(u\mu) = \Delta(\mu) &  \\
\partial_{t}(u\mu) + \partial_{x}(u^2\mu) + \partial_{x}(a\mu) = \Delta(u\mu) + \nu &  \\
\mu|_{t=0} = \mu_0,\;\; u\mu|_{t=0} = u_0\mu_0,\;\; a\mu|_{t=0} = 0,\;\; \nu|_{t=0} = \nu_0 &
\end{array}
\right.
\end{align}
Here, the unknowns $\mu$ and $\nu$ are fields of real Borel measures (with $\nu$ referred to as the force), and $u$, $a$ are real-valued functions. The fields $\mu$ and $u\mu$ are assumed to be càdlàg in time, with respective jumps $\Delta(\mu)$ and $\Delta(u\mu)$ (detailed in \ref{section gsp}). The initial data are given by $\mu_0$, $\nu_0$, and $u_0$, and the initial conditions are understood in the sense of weak convergence of measures as $t$ tends to zero.

Typically in the literature, the measure fields are assumed to be continuous, resulting in vanishing jumps. In that case, the system reduces to:
\begin{align}\label{equation gas}
\left\{
\begin{array}{ll}
\partial_{t}(\mu) + \partial_{x}(u\mu) = 0 &  \\
\partial_{t}(u\mu) + \partial_{x}(u^2\mu) + \partial_{x}(a\mu) = \nu &  \\
\mu|_{t=0} = \mu_0,\;\; u\mu|_{t=0} = u_0\mu_0,\;\; a\mu|_{t=0} = 0,\;\; \nu|_{t=0} = \nu_0 &
\end{array}
\right.
\end{align}

As discussed in
\cite{Anoverviewoncongestion,Priseencomptedelacongestion,Selforganizedhydrodynamics,NumericalsimulationsoftheEuler,Lesmilieuxgranulaires,Freecongestedtwophase,Amodelfortheevolution,Onafreeboundarybarotropic,MathematicalTopicsinFluidMechanics,Oncompactnessofsolutions,SingularlimitofaNavier-Stokes},
this system models the flow of a pressureless gas described by a mass density field $\mu$, an Eulerian velocity field $u$, and subject to an external force field $\nu$. A key feature of the model is the inclusion of a congestion pressure term $a \geq 0$, defined per unit mass. This term vanishes in regions of low density and becomes nonzero ($a \neq 0$) exclusively when the density reaches its maximal value, thereby capturing the onset of congestion effects.

In many works, this system is related to the sticky particle model, which describes particles moving freely along the real line at constant velocities in the absence of shocks. All shocks are completely inelastic and governed by the physical laws of conservation of mass and momentum. This model was introduced by Zeldovich \cite{A2,PeeblesThelargeScale} (\textcolor{blue}{1970-1980}) to explain the formation of megastructures in the universe, and was further developed by authors such as Kofman, Shandarin \cite{A19,A20}, and Carnevale, Pomeau, Young \cite{Statisticsofballistic} (\textcolor{blue}{1990}).

The case $\nu \equiv 0 \equiv a\mu$ has been extensively studied. Notable references include \cite{vpDermouneMoutsinga,Probalistic-interpretation-Dermoune,Probabilistic-interpretation-law-Dermoune,Convex-hullsMoutsinga,Generalized-variational-principles-E-Rykov-Sinai,LaxHyperbolicsystems,OelinikDiscontinuoussolutions,Burgers-Hopf,Hynd2019}, where $\mu$ represents the mass distribution field of the sticky dynamics, and $u$ is its Eulerian velocity field. 
In other notable studies, \cite{Aifang2022, Aifang2025} focus on the motion of a piston surrounded by sticky particles. Then $\mu$ represents the mass of particles accumulated on the piston and $u$ is its Eulerian velocity field.

For $\nu \not\equiv 0 \equiv a\mu$, interesting results have been obtained in \cite{A21,NzissilaForcedGas}. The solution proposed by Brenier et al. \cite{A21} (\textcolor{blue}{2013}) involves accelerating sticky particles between successive shocks. In this case, $\mu$ and $u$ are as above, and $\nu = \gamma\mu$, where $\gamma$ is the acceleration field. In contrast, Nzissila and Moutsinga \cite{NzissilaForcedGas} (\textcolor{blue}{2023}) retained the unaccelerated sticky particle model described earlier. They focused on the specific motion of particles they termed \textit{turbulent particles} (e.g., cluster endpoints), and \textit{redistributed} the total mass of the remaining particles over them. They obtained a weak solution where $\mu$ is the mass distribution field of the \textit{turbulent particle motion}, $u$ is its Eulerian velocity field, and $\nu$ measures the \textit{force of turbulence} (or shocks), concentrated at shock sites.

E et al. \cite{Generalized-variational-principles-E-Rykov-Sinai} (\textcolor{blue}{1996}) were the first to define the unaccelerated sticky particle model at the continuous level (diffuse $\mu_0$). Their key tool is the Generalized Variational Principle (GVP), which allows one to determine, at any time $t$, all clusters the initial locations of particles that have stuck together. Dermoune and Moutsinga \cite{Convex-hullsMoutsinga} (\textcolor{blue}{2008}) provided another construction of the unaccelerated sticky particle model at the continuum level. Their main tool is the use of appropriate convex hulls, which allow one to state the principal properties of the dynamics and, consequently, to solve the gas system and recover cluster characterization from the GVP.

The aim of this work is to formulate, at the discrete level, a new GVP for the accelerated sticky particle model. We also investigate the Dermoune-type ordinary differential equations associated with this model and their connection to the gas systems (\ref{equation gas disc}) and (\ref{equation gas}).

The rest of the paper is organized as follows.  
In section \ref{section asp}, we outline the main properties of the discrete model of accelerated sticky particles. Within this framework, in section \ref{section gvp}, and under the assumption that the initial acceleration is a non-increasing function of position, we generalize the work of \cite{Generalized-variational-principles-E-Rykov-Sinai} by establishing a generalized variational principle that characterizes clusters sets of initial positions of stuck particles at any time. Using well-defined functionals over these clusters, the entire dynamics of the system is reconstructed in section \ref{section Flow and Dermoune-type equation}. Here, positions, velocities, and accelerations are described by stochastic processes that solve ordinary differential equations of the Dermoune type (\ref{th Dermoune_gl}).  
Finally, section \ref{section gsp} is devoted to the gas system in its two distinct forms (\ref{equation gas disc}) and (\ref{equation gas}). We begin with the more familiar form (\ref{equation gas}), corresponding to continuous solutions. A weak solution is obtained following approaches established in the literature, based on the sticky particle mass distribution over positions and their velocity function. In this construction, $a \equiv 0$. In subsection \ref{section velocity_gsp}, the form (\ref{equation gas disc}), corresponding to càdlàg solutions, is addressed through a mass measure defined over particle velocities and through the acceleration function. This shift in perspective, i.e. formulating the dynamics in velocity space instead of position space, suppresses the force ($\nu \equiv 0$) and gives rise to a congestion-like term $a \geq 0$, which is in fact the conditional variance of acceleration given by velocity.

\tableofcontents

\section{Accelerated Sticky Particle Model}\label{section asp}

We consider $N$ particles located on the real line at positions $x_1 < x_2 < \cdots < x_N.$
Each particle is characterized by its initial position $x_j$, initial velocity $v_j$, initial mass $m_j$, and an initial acceleration $\theta_j$.

The particles evolve under uniform acceleration, meaning that each particle moves with constant acceleration until a collision occurs. Upon impact, particles undergo a soft collision and stick together, following the physical laws of conservation of mass, momentum, and force.

Consequently, any particle initially located at position $x_i$ belongs, at any time $t > 0$, to a cluster, a group of particles whose initial positions lie consecutively and have merged before or at time $t$. We denote this cluster by the set of initial positions $
x_{g(t)}, \dots, x_{d(t)}.$

Such a cluster behaves as a single composite particle with total mass:
\begin{equation}
	m_i(t) = \displaystyle\sum_{g(t)\leq j \leq d(t)}{m_j}\,\quad\forall i\in\{g(t),\dots,d(t)\}
\end{equation}
(due to the conservation law of mass), with the acceleration $\theta_i(t)$ such that
\begin{equation}
	m_i(t)\theta_i(t) = \displaystyle\sum_{g(t)\leq j \leq d(t)}{m_j\theta_j}
\end{equation}
(due to the conservation law of force), 
with velocity $v_i(t)$ such that
\begin{equation}
	m_i(t) v_i(t)= \displaystyle\sum_{g(t)\leq j \leq d(t)}{m_j\left(v_j + t\theta_j\right)}
\end{equation}
(obtained by the conservation of force and momentum).
And the new particle is at position $x_i(t)$ such that 
\begin{equation}
	m_i(t)x_i(t) = \displaystyle\sum_{g(t)\leq j \leq d(t)}{m_j\left(x_j + tv_j + \frac{1}{2}t^2\theta_j\right)}.
\end{equation}

\begin{proposition}\label{th dynamique}
\begin{description}
	\item(i) $\forall i \in\{g(t),\cdots,d(t)\} :$	
	\begin{align}
		v_i(t) &= \frac{\displaystyle\sum_{g(t)\leq j \leq d(t)}{m_j\left(v_j + t\theta_j\right)}}{\displaystyle\sum_{g(t)\leq j \leq d(t)}{m_j}}
		= \frac{\displaystyle\sum_{g(t)\leq j \leq d(t)}{m_jv_j }}{\displaystyle\sum_{g(t)\leq j \leq d(t)}{m_j}}+t\theta_i(t);\\
		x_i(t) &= \frac{\displaystyle\sum_{g(t)\leq j \leq d(t)}{m_j\left(x_j + tv_j + \frac{1}{2}t^2\theta_j\right)}}{\displaystyle\sum_{g(t)\leq j \leq d(t)}{m_j}}\\
		&= \frac{\displaystyle\sum_{g(t)\leq j \leq d(t)}{m_jx_j }}{\displaystyle\sum_{g(t)\leq j \leq d(t)}{m_j}}+tv_i(t)-\frac{1}{2}t^2\theta_i(t).
	\end{align}
	\item(2i) $\forall i \in\{1,\cdots,N\} :$	
	\begin{align}
		v_i(t) &= v_i+\int_{0}^{t}\theta_i(s)\partiel s;\\
		x_i(t) &= x_i+\int_{0}^{t}v_i(s)\partiel s.
	\end{align}
\end{description}
\end{proposition}

Note that any cluster can be represented as an interval of the form
$[\alpha, \beta] \cap \{x_1, x_2, \dots, x_N\}$,
where $\alpha = x_{g(t)}$ and $\beta = x_{d(t)}$. For simplicity, we will refer to such clusters as intervals $[\alpha, \beta]$ throughout the remainder of the paper.

\begin{Definition}
At each time $t\geq 0$, the partition of the $\{x_1,x_2,\dots,x_N\}$ into clusters formed before or at time $t$ is denoted by $\xi_t = \{G_1,\cdots,G_{N(t)}\}$.
\end{Definition}

We now interpret the dynamics of the system based on the following initial data:
\begin{itemize}
    \item the initial mass distribution, given by the discrete measure
    $P_0 = \sum_{j=1}^N m_j \delta_{x_j}$,
    where $m_j$ denotes the mass located at position $x_j$;
    
    \item the initial velocity function
    $u_0 : x_i \mapsto v_i$,
    
    assigning to each particle at position $x_i$ its initial velocity $v_i$;
    
    \item the initial acceleration function
    $\gamma_0 : x_i \mapsto \theta_i$,
    which associates to each particle its initial acceleration $\theta_i$.
\end{itemize}
\begin{Corollary}
	$\forall x_i\in[\alpha,\beta]\in\xi_t$,
	\begin{align}
		m_i(t)&=P_0([\alpha,\beta])\\
		\theta_i(t)&=\frac{\int_{[\alpha,\beta]}\gamma_{0}(x)\partiel P_0(x)}{P_0([\alpha,\beta])}\\
		v_i(t)&=\frac{\int_{[\alpha,\beta]}\left(u_0(x)+t\gamma_{0}(x)\right) \partiel P_0(x)}{P_0([\alpha,\beta])}\\
		x_i(t)&=\frac{\int_{[\alpha,\beta]}\left(x+tu_0(x)+\frac{t^2}{2}\gamma_{0}(x)\right) \partiel P_0(x)}{P_0([\alpha,\beta])}
	\end{align}
\end{Corollary}

This result shows that to construct the dynamics of accelerated sticky particles, it is enough to determine, at any time $t$, the set $\xi_t$ of clusters $[\alpha,\beta]$; we then apply the above formulae to obtain the mass, acceleration, velocity, and position of each cluster.

\section{Generalized Variational Principle}\label{section gvp}

In the case where the initial acceleration vanishes, i.e., $\gamma_0 \equiv 0$, a generalized variational principle (GVP), established in \cite{Generalized-variational-principles-E-Rykov-Sinai}, provides a constructive framework for determining the partition $\xi_t$ into clusters. This enables the construction of the sticky particle dynamics.

Our objective is to adapt this procedure to the more general case where $\gamma_0 \not\equiv 0$, i.e., when particles are subject to initial accelerations.

Let us begin by recalling this generalized variational principle (GVP). In \cite{Generalized-variational-principles-E-Rykov-Sinai}, any element $[\alpha,\beta]$ of $\xi_t$ is completely characterized by the following three properties :

\begin{description}
	\item[(i)] \begin{align*}
\forall y \in [\alpha,\beta] ,\hspace{1cm} \frac{\displaystyle\int_{[\alpha,y)}\displaystyle\left(x+tu_0(x)\displaystyle\right)dP_0(x)}{P_0([\alpha,y))} \geq \frac{\displaystyle\int_{(y,\beta]}\displaystyle\left(x+tu_0(x)\displaystyle\right)dP_0(x)}{P_0((y,\beta])}.
\end{align*}
\item[(ii)] \begin{align*} \forall y_1 < \alpha < y_2 :
\frac{\displaystyle\int_{[y_1,\alpha)}\displaystyle\left(x+tu_0(x)\displaystyle\right)dP_0(x)}{P_0([y_1,\alpha))} < \frac{\displaystyle\int_{[\alpha,y_2)}\displaystyle\left(x+tu_0(x)\displaystyle\right)dP_0(x)}{P_0([\alpha,y_2))} 
\end{align*}
\item[(iii)] \label{pbGVP}
\begin{align*} \forall y_1 < \beta < y_2 :
\frac{\displaystyle\int_{(y_1,\beta]}\displaystyle\left(x+tu_0(x)\displaystyle\right)dP_0(x)}{P_0((y_1,\beta])} < \frac{\displaystyle\int_{(\beta,y_2]}\displaystyle\left(x+tu_0(x)\displaystyle\right)dP_0(x)}{P_0((\beta,y_2])} \\
\end{align*}
\end{description}
Note that, in the case of a two-particle system where the particles collide at time $T$, the pair $\{x_1, x_2\}$ forms a cluster for all times $t > T$. Then property (iii) is equivalent to :
$$x_1+tv_1\geq x_2+tv_2,\quad\forall\,t\geq T\;.$$
If we adapt this property to the case where these particles are accelerated, we should have
\begin{align}\label{pbGVP accelere} 
	y_1(t):=x_1+tv_1+\frac{t^2}{2}\theta_1\geq x_2+tv_2+\frac{t^2}{2}\theta_2=:y_2(t),\quad\forall\,t\geq T\;.
\end{align}
We already know from the proposition \ref{th dynamique}, that
\begin{align*}
	y_1(t)=x_1(t)<x_2(t)=y_2(t)\,,&\quad\forall\,t<T\;;\\
	x_1(t)=x_2(t)=\frac{m_1y_1(t)+m_2y_2(t)}{m_1+m_2}\,,&\quad\forall\,t\geq T.
\end{align*}
To obtain the inequation (\ref{pbGVP accelere}), we quickly study the two parabolas $t\mapsto y_1(t), y_2(t)$; this leads to two interesting cases.
\begin{description}
	\item[Case 1 : $\mathbf{\theta_1<\theta_2}$.] The two parabolas meet at two times $t_1=T$ and $t_2>T$. Then we have
		\begin{align*}
		y_1(t)>y_2(t)\,,&\quad\forall\,t\in]T,t_2[\\
		y_1(t)<y_2(t)\,,&\quad\forall\,t>t_2.
		\end{align*}
The inequation (\ref{pbGVP accelere}) does not hold when $t>t_2$ (see figure \textbf{\textcolor{purple}{1}}).
	\item[Case 2 : $\mathbf{\theta_1>\theta_2}$.] The two parabolas meet at a single time (on $\mathbb{R}_+^*$) $t_1=T$. Then we have 
	\begin{align*}
		y_1(t)>y_2(t)\,,\quad\forall\,t>T.
	\end{align*} 
The inequation (\ref{pbGVP accelere}) holds for all $t>T$ (see figure \textbf{\textcolor{purple}{2}}).
\end{description}\begin{figure}[h]
	\centering
	\begin{minipage}[t]{0.45\textwidth}
		\includegraphics[width=\textwidth]{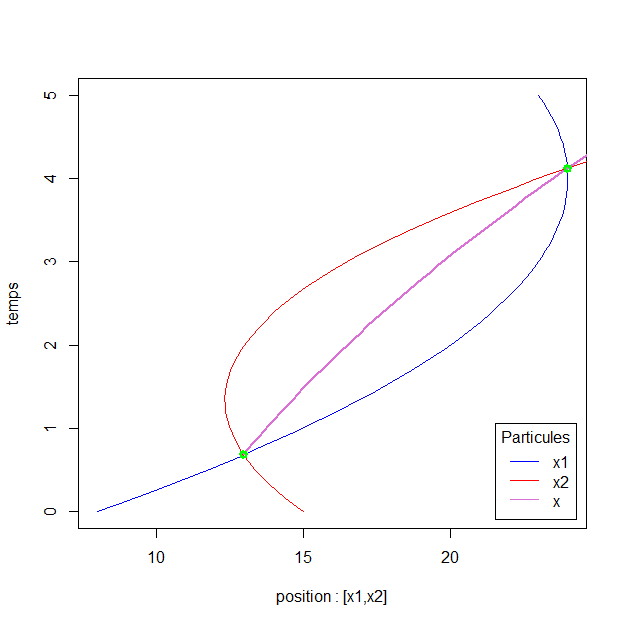}
	\caption{$\mathbf{\theta_1<\theta_2}$ : The inequation (\ref{pbGVP accelere}) does not hold when $t>t_2$.}
	\label{fig:Trajectoires1}
	\end{minipage}
\hfill
	\begin{minipage}[t]{0.45\textwidth}
		\includegraphics[width=\textwidth]{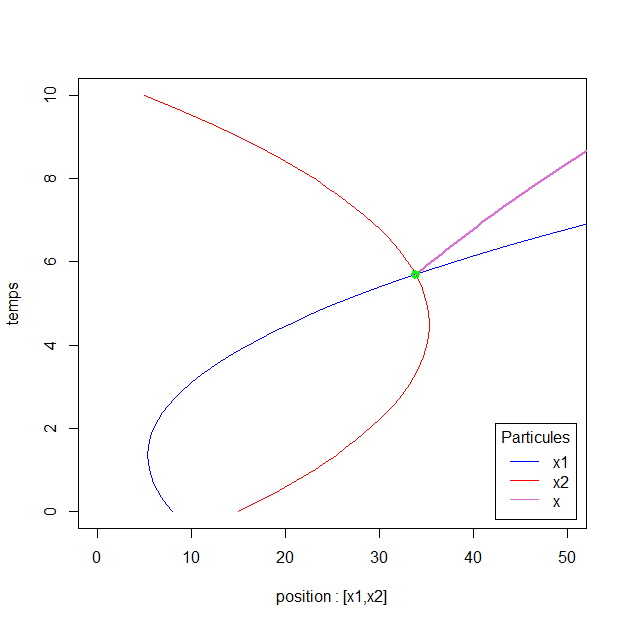}
	\caption{$\mathbf{\theta_1>\theta_2}$ : The inequation (\ref{pbGVP accelere}) is verified all the time after the shock.}
	\label{fig:Trajectoires2}
	\end{minipage}
\end{figure}

The preceding example compels us to discard the case $\theta_1 < \theta_2$. More generally, we exclude scenarios in which the function $\gamma_0$ is increasing, and henceforth restrict our attention to the case where $\gamma_0$ is non-increasing.

\begin{Theorem}\label{th gvpasp}
Assume that the function $\gamma_0$ is non-increasing. Let $t > 0$, and let $\xi_t$ denote the set of clusters at time $t$.
\begin{enumerate}
\item For all $y \in [\alpha, \beta] \in \xi_t$, the following inequality holds:
\begin{align}\label{eq gvpasp-choc}
\frac{\displaystyle\int_{[\alpha,y)}\left(\eta + t u_0(\eta) + \frac{t^2}{2} \gamma_0(\eta)\right) \, dP_0(\eta)}{P_0([\alpha,y))} 
\geq 
\frac{\displaystyle\int_{[y,\beta]}\left(\eta + t u_0(\eta) + \frac{t^2}{2} \gamma_0(\eta)\right) \, dP_0(\eta)}{P_0([y,\beta])}
\end{align}

\item \textbf{Left endpoint condition} $(GVP_{\gamma_g})$: $\alpha$ is the left endpoint of an element of $\xi_t$ if and only if, for all $y_1 < \alpha < y_2$, we have:
\begin{align}
\frac{\displaystyle\int_{[y_1,\alpha)}\left(\eta + t u_0(\eta) + \frac{t^2}{2} \gamma_0(\eta)\right) \, dP_0(\eta)}{P_0([y_1,\alpha))} 
< 
\frac{\displaystyle\int_{[\alpha,y_2)}\left(\eta + t u_0(\eta) + \frac{t^2}{2} \gamma_0(\eta)\right) \, dP_0(\eta)}{P_0([\alpha,y_2))}
\end{align}

\item \textbf{Right endpoint condition} $(GVP_{\gamma_d})$: $\beta$ is the right endpoint of an element of $\xi_t$ if and only if, for all $y_1 < \beta < y_2$, we have:
\begin{align}
\frac{\displaystyle\int_{(y_1,\beta]}\left(\eta + t u_0(\eta) + \frac{t^2}{2} \gamma_0(\eta)\right) \, dP_0(\eta)}{P_0((y_1,\beta])} 
< 
\frac{\displaystyle\int_{(\beta,y_2]}\left(\eta + t u_0(\eta) + \frac{t^2}{2} \gamma_0(\eta)\right) \, dP_0(\eta)}{P_0((\beta,y_2])}
\end{align}
\end{enumerate}
\end{Theorem}

\begin{Remark}
For any time $t > 0$, the elements of $\xi_t$ are precisely the intervals $[\alpha, \beta]$ such that $\alpha$ satisfies $(GVP_{\gamma_g})$ at time $t$, $\beta$ satisfies $(GVP_{\gamma_d})$ at time $t$, and no point in $(\alpha, \beta)$ satisfies either $(GVP_{\gamma_g})$ or $(GVP_{\gamma_d})$.
\end{Remark}
We thus recover the classical variational principle of \cite{Generalized-variational-principles-E-Rykov-Sinai} in the particular case where $\gamma_0 \equiv 0$. Our result therefore constitutes a true generalization of the framework developed in \cite{Generalized-variational-principles-E-Rykov-Sinai}. 

To establish \ref{th gvpasp}, we first require the following elementary lemma:

\begin{Lemma}
Let $Q_1,Q_2$ be two quadratic functions such that $Q_1\not\equiv Q_2$ and $Q_1''\geq Q_2''$. If there exists $t_0<t_1$ such that $Q_1(t_0)\leq Q_2(t_0)$ and $Q_1(t_1)\geq Q_2(t_1)$, then $Q_1(s)>Q_2(s)$, $\forall,s> t_1$.
\end{Lemma}

\noindent\textbf{Proof of the theorem \ref{th gvpasp} :}
\begin{enumerate}
 \item We proceed to a demonstration by recurrence on the successive shock times $t_{1}(y)< t_{2}(y)< \cdots< t_{n}(y)\leq t$ of the particle $y$. We simply denote $t_k$ to simplify the notation. At each shock time $t_k$, $y$ is in a cluster $[\alpha_k,\beta_k]$ ($k=1,\dots,n$). We have $[\alpha,\beta]=[\alpha_n,\beta_n]$. Consider the quadratic functions $Q_k$ and $R_k$ defined by
\begin{align*}
	Q_{k}(s) &= \displaystyle\frac{\displaystyle\int_{[\alpha_k,y)}\eta dP_0(\eta)}{P_0([\alpha_k,y))}
	+s\displaystyle\frac{\displaystyle\int_{[\alpha_k,y)}u_{0}(\eta)dP_0(\eta)}{P_0([\alpha_k,y))} +\frac{s^2}{2}\displaystyle\frac{\displaystyle\int_{[\alpha_k,y)}\gamma_{0}(\eta)dP_0(\eta)}{P_0([\alpha_k,y))}\;,\\
	R_{k}(s) &= \displaystyle\frac{\displaystyle\int_{[y,\beta_k]}\eta dP_0(\eta)}{P_0([y,\beta_k])}
	+s\displaystyle\frac{\displaystyle\int_{[y,\beta_k]}u_{0}(\eta)dP_0(\eta)}{P_0([y,\beta_k])} +\frac{s^2}{2}\displaystyle\frac{\displaystyle\int_{[y,\beta_k]}\gamma_{0}(\eta)dP_0(\eta)}{P_0([y,\beta_k])}\;.
\end{align*}
For $k=1$, forall $x_i\in[\alpha_1,\beta_1]$, we have  
$$x_i+t_1u_{0}(x_i)+\frac{t_{1}^{2}}{2}\gamma_{0}(x_i) = y+t_1u_{0}(y)+\frac{t_{1}^{2}}{2}\gamma_{0}(y).$$ 
Then
\begin{align*}
	\displaystyle\frac{\displaystyle\int_{[\alpha,y)}\left(\eta+t_1u_0(\eta)+\frac{t_1^2}{2}\gamma_0(\eta)\right)dP_0(\eta)}{P_0([\alpha,y))} = \displaystyle\frac{\displaystyle\int_{[y,\beta]}\left(\eta+t_1u_0(\eta)+\frac{t_1^2}{2}\gamma_0(\eta)\right)dP_0(\eta)}{P_0([y,\beta])}\;,
\end{align*} 
which means that $Q_{1}(t_1)=R_{1}(t_1)$. But we have $Q_{1}(0)<y\leq R_{1}(0)$.
Since the function $\gamma_0$ is decreasing, we have $Q_{1}^{''}\geq R_{1}^{''}$. Then, according to the previous lemma, we have $Q_{1}(s)>R_{1}(s)$, $\forall\,s> t_1$. Thus, as long as $[\alpha_1,\beta_1]$ remains a cluster, the relation (\ref{eq gvpasp-choc}) remains true for $[\alpha_1,\beta_1]$.
Now suppose that for $k\geq1$, we have $Q_{k}(s)>R_{k}(s)$, $\forall\,s> t_k$. Let's show that $Q_{k+1}(s)>R_{k+1}(s)$, $\forall\,s> t_{k+1}$.
Recall that at time $t_k$ of shock, $y$ is in cluster $[\alpha_{k},\beta_{k}]$. At time $t_{k+1}$, there are $g$ consecutive clusters to the left of $y$ and $d$ consecutive clusters to the right of $y$ which collide with $[\alpha_{k},\beta_{k}]$ (for any integers $g,d$). Let's note these clusters $[\alpha_{1}^{'},\beta_{1}^{'}]$, $[\alpha_{2}^{'},\beta_{2}^{'}]$, $\cdots$, $[\alpha_{g}^{'},\beta_{g}^{'}]$ and
$[\alpha_{1}^{''},\beta_{1}^{''}]$, $[\alpha_{2}^{''},\beta_{2}^{''}]$, $\cdots$, $[\alpha_{d}^{''},\beta_{d}^{''}]$. We have inevitably$\alpha_{1}'=\alpha_{k+1}$ et $\beta_{d}^{''}=\beta_{k+1}$ (see Figures \textbf{\textcolor{purple}{3}} and \textbf{\textcolor{purple}{4}}). 
Then $\eta(s)=\eta+su_0(\eta)+\frac{s^2}{2}\gamma_0(\eta)$, we have
$\forall i,j$ :
\begin{align*}
\frac{\displaystyle\int_{[\alpha_{i}^{'},\beta_{i}^{'}]}\eta(t_{k+1})dP_0(\eta)}{P_0([\alpha_{i}^{'},\beta_{i}^{'}])} = \frac{\displaystyle\int_{[\alpha_{j}^{''},\beta_{j}^{''}]}\eta(t_{k+1})dP_0(\eta)}{P_0([\alpha_{j}^{''},\beta_{j}^{''}])} = \frac{\displaystyle\int_{[\alpha_{k},\beta_{k}]}\eta(t_{k+1})dP_0(\eta)}{P_0([\alpha_{k},\beta_{k}])}=:z\;,\\
z=\frac{\displaystyle\int_{[\alpha_{k+1},\beta_{k+1}]}\eta(t_{k+1})dP_0(\eta)}{P_0([\alpha_{k+1},\beta_{k+1}])} = \frac{\displaystyle\int_{[\alpha_{k+1},\alpha_{k})}\eta(t_{k+1})dP_0(\eta)}{P_0([\alpha_{k+1},\alpha_{k}))}
\leq \frac{\displaystyle\int_{[\alpha_{k},y)}\eta(t_{k+1})dP_0(\eta)}{P_0([\alpha_{k},y))}\;.
\end{align*}

The last two equalities are obtained by calculating the barycentre of several points equal to $z$. The inequality is obtained using the recurrence hypothesis. 
The mean of the last two terms is then obtained:
\begin{align*}
	\frac{\displaystyle\int_{[\alpha_{k+1},y)}\eta(t_{k+1})dP_0(\eta)}{P_0([\alpha_{k+1},y))}\geq 
	\frac{\displaystyle\int_{[\alpha_{k+1},\beta_{k+1}]}\eta(t_{k+1})dP_0(\eta)}{P_0([\alpha_{k+1},\beta_{k+1}])} .
\end{align*}
Since $y\in[\alpha_{k+1},\beta_{k+1}]$, we inevitably have
\begin{align*}
\frac{\displaystyle\int_{[y,\beta_{k+1}]}\eta(t_{k+1})dP_0(\eta)}{P_0([y,\beta_{k+1}])}\leq z
\leq \frac{\displaystyle\int_{[\alpha_{k+1},y)}\eta(t_{k+1})dP_0(\eta)}{P_0([\alpha_{k+1},y))}\;;
\end{align*}
which means that $Q_{k+1}(t_{k+1})\geq R_{k+1}(t_{k+1})$. The argument used above (using the previous lemma) shows that $Q_{k+1}(s)>R_{k+1}(s)$ $\forall\,s> t_{k+1}$.
The recurrence formula is therefore true for all $k$. So, as long as $[\alpha_n,\beta_n]$ remains a cluster, the relation (\ref{eq gvpasp-choc}) remains true for $[\alpha_n,\beta_n]$. This gives (\ref{eq gvpasp-choc}) for $[\alpha,\beta]$ and $t$.\\
It is this last property that will be used for the rest of the proof, in addition to simple calculations of barycentres.

\begin{figure}[h]
	\centering
 		\includegraphics[height=4cm]{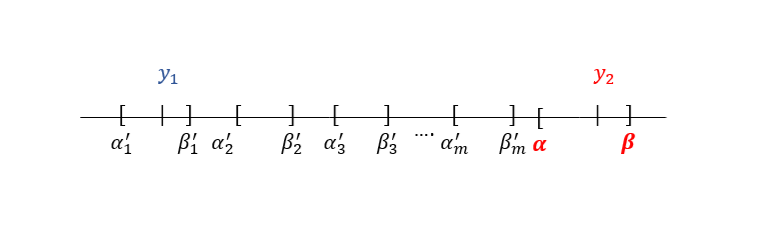}
	\caption{Consecutive clusters on the left of $\alpha$}\label{fig:Extrémité_gauche}

\hfill
		\includegraphics[height=4.2cm]{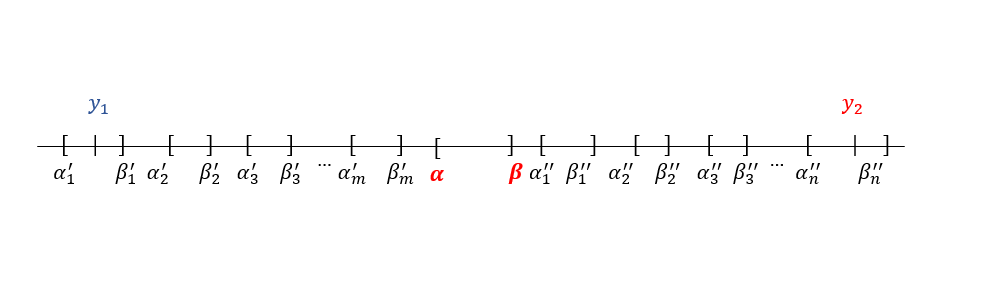}
	\caption{Consecutive clusters around $\alpha$}\label{fig:Extrémité_droite}
\end{figure}
%
\item Consider $y_1<\alpha<y_2$, with $[\alpha,\beta]\in\xi_t$. Since $y_1<\alpha$, there are consecutive clusters (elements of $\xi_t$) on the left of $\alpha$ which we denote
$[\alpha_{1}',\beta_{1}^{'}]$, $[\alpha_{2}^{'},\beta_{2}^{'}]$, $\cdots$, $[\alpha_{g}^{'},\beta_{g}^{'}]$ such that 
$y_1\in[\alpha_{1}',\beta_{1}^{'}]$ (see Figure \ref{fig:Extrémité_gauche}). For all $1<j\leq g$, we have the positions at time $t$ :
\begin{align}\label{etoile1}
\displaystyle\frac{\displaystyle\int_{[\alpha_{1}^{'},\beta_{1}^{'}]}\eta(t) dP_0(\eta)}{P_0([\alpha_{1}^{'},\beta_{1}^{'}])}<\displaystyle\frac{\displaystyle\int_{[\alpha_{j}^{'},\beta_{j}^{'}]}\eta(t) dP_0(\eta)}{P_0([\alpha_{j}^{'},\beta_{j}^{'}])}<\displaystyle\frac{\displaystyle\int_{[\alpha,\beta]}\eta(t) dP_0(\eta)}{P_0([\alpha,\beta])}\;.
\end{align}
Since $y_1\in [\alpha_{1}^{'},\beta_{1}^{'}]$, we apply the assertion 1 and obtain :
\begin{align*}
\displaystyle\frac{\displaystyle\int_{[y_{1},\beta_{1}^{'}]}\eta(t)dP_0(\eta)}{P_0([y_1,\beta_{1}^{'}])}\leq\displaystyle\frac{\displaystyle\int_{[\alpha_{1},y_{1})}\eta(t)dP_0(\eta)}{P_0([\alpha_{1}^{'},y_{1}))}\;;
\end{align*}
which leads to
\begin{align}\label{etoile2}
\displaystyle\frac{\displaystyle\int_{[y_{1},\beta_{1}^{'}]}\eta(t)dP_0(\eta)}{P_0([y_1,\beta_{1}^{'}])}\leq
\displaystyle\frac{\displaystyle\int_{[\alpha_{1}^{'},\beta_{1}^{'}]}\eta(t) dP_0(\eta)}{P_0([\alpha_{1}^{'},\beta_{1}^{'}])}
\leq
\displaystyle\frac{\displaystyle\int_{[\alpha_{1}^{'},y_1)}\eta(t) dP_0(\eta)}{P_0([\alpha_{1}^{'},y_1))}
\;.
\end{align}
Therefore, from (\ref{etoile1}) and (\ref{etoile2}), we get for all $1<j\leq g$ :
\begin{align*}
\displaystyle\frac{\displaystyle\int_{[y_{1},\beta_{1}^{'}]}\eta(t)dP_0(\eta)}{P_0([y_1,\beta_{1}^{'}])}<
\displaystyle\frac{\displaystyle\int_{[\alpha_{j}^{'},\beta_{j}^{'}]}\eta(t) dP_0(\eta)}{P_0([\alpha_{j}^{'},\beta_{j}^{'}])}<
\displaystyle\frac{\displaystyle\int_{[\alpha,\beta]}\eta(t) dP_0(\eta)}{P_0([\alpha,\beta])}
\end{align*}
A calculation of the barycentre shows that
\begin{align}\label{etoile3}
	\displaystyle\frac{\displaystyle\int_{[y_{1},\alpha[}\eta(t)dP_0(\eta)}{P_0([y_1,\alpha[)}<
	\displaystyle\frac{\displaystyle\int_{[\alpha,\beta]}\eta(t) dP_0(\eta)}{P_0([\alpha,\beta])}
\end{align}
Now consider $y_2>\alpha$. There are two cases.
\begin{description} 
	\item[If $\mathbf{y_2\leq\beta}$] (see Figure \ref{fig:Extrémité_gauche}), then in the same way as (\ref{etoile2}) : 
\begin{align*}
\displaystyle\frac{\displaystyle\int_{[y_{2},\beta]}\eta(t)dP_0(\eta)}{P_0([y_2,\beta])}\leq
\displaystyle\frac{\displaystyle\int_{[\alpha,\beta]}\eta(t) dP_0(\eta)}{P_0([\alpha,\beta])}
\leq
\displaystyle\frac{\displaystyle\int_{[\alpha,y_2)}\eta(t) dP_0(\eta)}{P_0([\alpha,y_2))}
\;.
\end{align*}
The last inequality, combined with (\ref{etoile3}), gives the required result.
\item[If $\mathbf{y_2>\beta}$] (see Figure \ref{fig:Extrémité_droite}), 
there are consecutive clusters (elements of $\xi_t$) on the right of $\beta$ which we denote $[\alpha_{1}^{''},\beta_{1}^{''}]$, $[\alpha_{2}^{''},\beta_{2}^{''}]$, $\cdots$, $[\alpha_{d}^{''},\beta_{d}^{''}]$ such that $y_2\in[\alpha_{d}^{''},\beta_{d}^{''}]$. Then
\begin{align*}
	\displaystyle\frac{\displaystyle\int_{[\alpha_{d}^{''},\beta_{d}^{''}]}\eta(t) dP_0(\eta)}{P_0([\alpha_{d}^{''},\beta_{d}^{''}])}
	\leq
	\displaystyle\frac{\displaystyle\int_{[\alpha_{d}^{''},y_2)}\eta(t) dP_0(\eta)}{P_0([\alpha_{d}^{''},y_2))}
	\;.
\end{align*}

Barycenter arguments based on the successive positions of clusters then lead to\begin{align*}
	\displaystyle\frac{\displaystyle\int_{[\alpha,\beta]}\eta(t) dP_0(\eta)}{P_0([\alpha,\beta])}
	\leq
	\displaystyle\frac{\displaystyle\int_{[\alpha_{1}^{''},y_2)}\eta(t) dP_0(\eta)}{P_0([\alpha_{1}^{''},y_2))}
	\;;
\end{align*}
which, combined with (\ref{etoile3}), gives the desired result.
\end{description} 

\item Similar to the previous case.
\end{enumerate}
\hfill$\square$

\section{Flow and Equation of the Dermoune Type}\label{section Flow and Dermoune-type equation}

Positions, velocities, and accelerations can be defined as functions over clusters $[\alpha, \beta]$ at time $t$ as follows. For all $x_i \in [\alpha, \beta]$, we set:
\begin{align*}
\phi(x_i, t) &= x_i(t) = \frac{\displaystyle\int_{[\alpha,\beta]} \left[\eta + t u_0(\eta) + \frac{1}{2} t^2 \gamma_0(\eta)\right] \, dP_0(\eta)}{P_0([\alpha,\beta])}, \\
v(x_i, t) &= v_i(t) = \frac{\displaystyle\int_{[\alpha,\beta]} \left[u_0(\eta) + t \gamma_0(\eta)\right] \, dP_0(\eta)}{P_0([\alpha,\beta])}, \\
\theta(x_i, t) &= \theta_i(t) = \frac{\displaystyle\int_{[\alpha,\beta]} \gamma_0(\eta) \, dP_0(\eta)}{P_0([\alpha,\beta])}.
\end{align*}

The mass distribution over particle positions at time $t$ is given by the discrete measure $P_t := \sum m_j \delta_{x_j(t)}.$
Note that if $x_i(t) = x_j(t)$, then $v(x_i, t) = v(x_j, t)$ and $\theta(x_i, t) = \theta(x_j, t)$. This implies the existence of functions $u$ and $\gamma$ such that:
\begin{align*}
v(x_i, t) &= u(x_i(t), t), \\
\theta(x_i, t) &= \gamma(x_i(t), t).
\end{align*}

As in \cite{Convex-hullsMoutsinga}, and by construction, the properties established in \ref{section asp} translate into the following form :

\begin{proposition}[Generalized Dermoune equation]\label{th Dermoune_gl}
Let $X_0$ be a random variable with distribution $P_0$ and let the process $X$ be defined by $X_t = \phi(X_{0},t)$. 
The process $X$ is absolutely continuous and the processes $t\mapsto u(X_t,t),\,\gamma(X_t,t)$ are càdlàg in $\R_+^*$. They all three processes are right-continuous at 0. 
Moreover, the processes $t\mapsto X_t$, $u(X_t,t)$ have right-hand derivatives given by
\begin{align}\label{Eq3}
	\frac{\partiel X_t}{\partiel t} = \displaystyle\mathbb{E}\left[u_0(X_0) + t\gamma_0(X_0)\Big|X_t\right]\,.
\end{align}
\begin{align}\label{Eq5}
	\frac{\partiel u(X_t,t)} {\partiel t}= \displaystyle\mathbb{E}\left[\gamma_0(X_0)|X_t\right]\,.
\end{align}
In addition,  for all $t$ :
\begin{align}
	X_t= \displaystyle\mathbb{E}\left[X_0 + tu_0(X_0) + \frac{t^2}{2}\gamma_0(X_0)\Big|X_t\right]\,,
\end{align}
\begin{align}\label{Eq4}
	u(X_t,t) = \displaystyle\mathbb{E}\left[u_0(X_0) + t \gamma_0(X_0)|X_t\right]\;,
\end{align}
\begin{align}\label{Eq6}
	\gamma(X_t,t) = \displaystyle\mathbb{E}\left[\gamma_0(X_0)|X_t\right]\,.
\end{align}
\end{proposition}
\begin{Remark}\label{rem_cle Dermoune}
In the particular case where $\gamma_0 \equiv 0$, Dermoune derives (\ref{Eq3}) as presented in \cite{Probabilistic-interpretation-law-Dermoune, Probalistic-interpretation-Dermoune}. For this reason, we refer to the general form (\ref{Eq3}) as the \emph{generalized Dermoune equation}.

Moreover, when considering the simplified form (\ref{Eq4}) under the assumption $\gamma_0 \equiv 0$, Dermoune provides a probabilistic solution to the pressureless gas system (\ref{Eq2}). This solution is constructed via the velocity field $u$, the density $\rho(\cdot, t) = \text{Law}(X_t)$, and the vanishing acceleration field $\gamma \equiv 0$.
\end{Remark}

Note also that for the unaccelerated sticky particle model, equation (\ref{Eq3}) was also highly useful in \cite{Hynd2019} to get the pressureless gas system (\ref{Eq2}). 

\section{Pressure less accelerated gas}\label{section gsp}  

We now consider a model of gas particles without pressure. The gas goes to the liquefaction phase according to the accelerated sticky particle model described above. We show the link with the pressureless gas system of equations (\ref{equation gas disc}). All the the measures are supposed finite on compact subsets of $\R$. The derivatives of measure fields are understood is the sense of distributions. We need the limits, in the sense of weak convergence, for any field $m$ of measures :
\begin{align*}
	m(\cdot,t^+)=\lim_{\underset{s<t}{s\to t}}m(\cdot,t)\,,\quad m(\cdot,t^-)=\lim_{\underset{s>t}{s\to t}}m(\cdot,t)\,.
\end{align*}
Thus, we have new fields $m^+:=(m(\cdot,t^+),\,t>0)$ and $m^-:=(m(\cdot,t^-),\,t>0)$.
The field $m$ is said càdlàg if the above limit fields exist and $m^+=m$. We can define the signed measure on $\R\times\R_+^*$, related to the jumps of the càdlàg field $m$ :
\begin{align*}
	\Delta (m)=\sum_{t\in\R_+^*}(m^+-m^-)(\cdot,t)\delta_t\;,
\end{align*}
with $\delta_t$ the Dirac measure at $t$. It is a  well defined sum of product measures. Indeed, for all continuous function $f$ with compact support and all $0<t_1<t_2 : $
\begin{align*}
	\int_{]t_1,t_2]}\int_{\R}f(x)\Delta (m)(\partiel x,\partiel t)=
	\sum_{t_1< t\leq t_2}\int f(x)[m(\partiel x,t^+)-m(\partiel x,t^-)]\;.
\end{align*}
The function $t\mapsto \int f(x)m(\partiel x,t)$ is càdlàg with left hand limits $\int f(x)m(\partiel x,t^-)$, for all $t$. Then, there exits an at-most countable subset $D$ (depending on $f$) such that 
\begin{align*}
	\int_{]t_1,t_2]}\int_{\R}f(x)\Delta (m)(\partiel x,\partiel t)=
	\sum_{t\in ]t_1,t_2]\cap D}\int f(x)[m(\partiel x,t^+)-m(\partiel x,t^-)]\;.
\end{align*}

Now we can give the definition of a weak solution of (\ref{equation gas disc}). We use the space  
$C_{0}^{1}(\mathbb{R})$ of real functions which have continuous derivatives and compact support.

\begin{Definition}
	The family $(\mu,u\mu,a\mu,\nu)$ of Borel measure fields is said to be a weak solution of (\ref{equation gas disc}) if it is weakly càdlàg with respect to time and if for all functions $f \in C_{0}^{1}(\mathbb{R})$ and all $0<t_1<t_2 : $
	\begin{align*}
		(E1)\quad	\displaystyle\int_{\R} &f(x)\mu(\partiel x,t_2) - \displaystyle\int_{\R} f(x)\mu(\partiel x,t_1) = \displaystyle\int_{t_1}^{t_2}\int_{\R} f'(x)u(x,t)\mu(\partiel x,t) \partiel t\\
			&+\int_{]t_1,t_2]}\int_{\R}f(x)\Delta (\mu)(\partiel x,\partiel t)\,,\\
		(E2)\quad	\displaystyle\int_{\R} &f(x)u(x,t_2)\mu(\partiel x,t_2) - 
		\displaystyle\int_{\R} f(x)u(x,t_1)\mu(\partiel x,t_1)\\
		= &\int_{t_1}^{t_2}\int_{\R} f'(x)[u^2(x,t)+a(x,t)]\mu(\partiel x,t)\partiel t
		+\int_{t_1}^{t_2}\int_{\R} f(x)\nu(\partiel x,t)\partiel t\\
		&+\int_{]t_1,t_2]}\int_{\R}f(x)\Delta (u\mu)(\partiel x,\partiel t)\,.
	\end{align*}
	Moreover, for all continuous $g$ with compact support and as $t$ tends to zero :
	\begin{align*}
		(E3)\quad \quad\quad \quad  	\displaystyle\int g(x)\mu(\partiel x,t) &\rightarrow \displaystyle\int g(x)\mu_{0}(\partiel x)\\
		(E4)\quad\;
		\displaystyle\int g(x)a(x,t)\mu(\partiel x,t) &\rightarrow \displaystyle\int g(x)a_{0}(x)\mu_{0}(\partiel x)\\
		(E5)\quad\;
		 \displaystyle\int g(x)u(x,t)\mu(\partiel x,t) &\rightarrow \displaystyle\int g(x)u_{0}(x)\mu_{0}(\partiel x)\\
		(E6)\quad\quad \quad \quad  \displaystyle\int g(x)\nu(\partiel x,t) &\rightarrow \displaystyle\int g(x)\nu_{0}(\partiel x)\,.
	\end{align*}
\end{Definition}

\subsection{Forced pressureless gas system without congestion-like effect} 
We begin in the setting of continuous solutions of (\ref{equation gas disc}). The solution is given by the particle mass distribution $\rho$ over positions, the force $\nu$ is absolutely continuous with respect to $\rho$ and the is no congestion-like effect. Precisely, we have to solve the simplified form, on $\R\times\R_+^*$ :
\begin{align}\label{Eq2}
	\begin{cases}
		\partial_t(\rho) + \partial_x(u\rho) = 0\\
		\partial_t(u\rho) + \partial_x(u^2\rho) = \gamma\rho\\
		\rho_t|_{t=0}=\rho_0\,,\quad u\rho|_{t=0}=u_0\rho_0\,,\quad\gamma\rho|_{t=0}=\gamma_0\rho_0\,,
	\end{cases}
\end{align}
for given probability $\rho_0$ and Borel functions $u_0$ and $\gamma_0$. Not that such a link has already been published in the literature (see for example \cite{A21}).
The following result shows that the solution can be obtained using probabilistic tools, by solving a Dermoune-type equation.

\begin{Lemma}[Probabilistic resolution]
	If there are regular fields of Borel functions $u,\gamma$ and a real stochastic process $X$ solution of
	\begin{align*}
	dX_t = u(X_t,t)d&t=\mathbb{E}[u_0(X_0)+t\gamma_0(X_0)|X_t]dt, \hspace{0.2cm} \law(X_0) = \rho_0,\\
	\gamma(X_t,t)&=\mathbb{E}[\gamma_0(X_0)|X_t],
\end{align*}
then a weak solution of (\ref{Eq2}) is obtained by ($\rho(.,t)=\mathcal{L}(X_t)$, $u\rho,\gamma\rho$).
\end{Lemma}
We can now give the solution. We consider a discrete sticky particle dynamics of initial mass distribution $\rho_0$, initial velocity function $u_0$ and initial acceleration function $\gamma_0$. Consider the trajectories defined by $t\mapsto\phi(x,t)$ (in section \ref{section Flow and Dermoune-type equation}).
\begin{Theorem}
If $X_0$ is a real random variable of law $\rho_0$, then the process
$X$ defined by $X_t:=\phi(X_0,t)$ solves the problem.
\end{Theorem}
Indeed, Proposition \ref{th Dermoune_gl} shows that this process satisfies the previous lemma.

\noindent\textbf{Proof of the lemma.}
	For all $f\in C_{0}^{1}(\mathbb{R})$ and $0<t_1<t_2,$ we have
	\begin{align*}
		\int f(x)\rho(dx,t_2) &- \int f(x)\rho(dx,t_1) = \mathbb{E}[f(X_{t_2})-f(X_{t_1})]\\
		&= \mathbb{E}\left[\int_{t_1}^{t_2}f'(X_t)dX_t\right]\\
		&= \mathbb{E}\left[\int_{t_1}^{t_2}f'(X_t)u(X_t,t)dt\right]\\
		&= \int_{t_1}^{t_2}\mathbb{E}\left[f'(X_t)u(X_t,t)\right]dt\\
		&= \int_{t_1}^{t_2}\int f'(x)u(x,t)d\rho(x,t)dt
	\end{align*}
	This concludes the first point. Now, for the second point, we define
	\begin{eqnarray*}
		\delta=\int f(x)u(x,t_2)d\rho(x,t_2)dt-\int f(x)u(x,t_1)d\rho(x,t_1)dt.
	\end{eqnarray*}
	We have
	\begin{align*}
		\delta&=\mathbb{E}[f(X_{t_2})u(X_{t_2},t_2)]-\mathbb{E}[f(X_{t_1})u(X_{t_1},t_1)]\\
		&=\mathbb{E}[f(X_{t_2})(u_0(X_0)+t_2\gamma_0(X_0))-f(X_{t_1})(u_0(X_0)+t_1\gamma_0(X_0))]\\
		&=\mathbb{E}[g(t_2)-g(t_1)]=\mathbb{E}\left[\int_{t_1}^{t_2}\partiel g(t)\right]	\;.
	\end{align*}
Define $g(t)=f(X_{t})(u_0(X_0)+t\gamma_0(X_0))$ and compute 
\begin{align*}
\partiel g(t)=(u_0(X_0)+t\gamma_0(X_0))f'(X_t)\partiel X_t+f(X_t)\gamma_0(X_0)\partiel t&\\
=(u_0(X_0)+t\gamma_0(X_0))f'(X_t)u(X_t,t)\partiel t+f(X_t)\gamma_0(X_0)\partiel t\;.&
\end{align*}
Then
	\begin{align*}
		\delta
		=&\int_{t_1}^{t_2}\mathbb{E}[f'(X_t)(u_0(X_0)+t\gamma_0(X_0))u(X_t,t)]dt\\
		&+\int_{t_1}^{t_2}\mathbb{E}[f(X_t)\gamma_0(X_0)]dt\\
		=&\int_{t_1}^{t_2}\mathbb{E}[f'(X_t)\mathbb{E}[(u_0(X_0)+t\gamma_0(X_0))|X_t]u(X_t,t)]dt\\
		&+\int_{t_1}^{t_2}\mathbb{E}[f(X_t)\mathbb{E}[\gamma_0(X_0)|X_t]]dt\\
		=&\int_{t_1}^{t_2}\mathbb{E}[f'(X_t)u(X_t,t)u(X_t,t)]dt\\
		&+\int_{t_1}^{t_2}\mathbb{E}[f(X_t)\gamma(X_t,t)]dt\\
		=&\int_{t_1}^{t_2}\int f'(x)u^{2}(x,t)\rho(dx,t)dt\\
		&+\int_{t_1}^{t_2}\int f(x)\gamma(x,t)d\rho(x,t)dt.
	\end{align*}

The third point concerns the weak convergence of the laws. It follows from the right-continuity of the process $X$ at 0.

The fourth point is less straightforward and concerns the weak convergence of $u\rho$ towards $u_0\rho_0$ as $t\to0$. 
For any bounded and continuous function $g$, the process $t\mapsto g(X_{t})u(X_{t},t)$ is bounded and is right-continuous at 0. So we can use the dominated convergence theorem :
	\begin{align*}
	\lim_{t\to 0^{+}} \int g(x)u(x,t)d\rho(x,t) &=\lim_{t\to 0^{+}} \mathbb{E}[g(X_{t})u(X_{t},t)]
	= \mathbb{E}[g(X_{0})u(X_{0},0)]\\
	&= \mathbb{E}[g(X_{0})u_0(X_{0})]
	= \int g(x)u_{0}(x)d\rho_{0}(x)\,.
	\end{align*}
The same technique is used for the last convergence.		\hfill$\square$\\
	
\subsection{Unforced pressureless gas system with congestion-like effect}\label{section velocity_gsp}
We end this paper with an important solution of gas system from the velocity process $t\mapsto V_t:=u(X_t,t)$. The key is the  equation of $V$, stated in Proposition \ref{th Dermoune_gl}, which is close to Dermoune's. This motivates an adaptation of the proof of Proposition 3 to the velocity process. However, we face a major difficulty : the fact that $V$ is not continuous. 
This is why we get the new form (\ref{equation gas disc}) of gas equations, in the settings of càdlàg fields of measures. Our solution is given by the particle mass distribution over velocities. Then, a congestion-like term merges, due to the fact that acceleration may vary, for a given velocity. Moreover, and rather unexpectedly, the force turns out to be zero (in this space of velocities). 

To achieve our goal, we consider the processes defined in section \ref{section Flow and Dermoune-type equation} and define
\begin{eqnarray*}
	\mu(\cdot,t)=\law(V_t)\,,&\; \mu^-(\cdot,t)=\law(V_{t^-})\,,&\; \mu_0=\law(V_0)\,,\\
	w(x,t)=\mathbb{E}[\Gamma_{0}|V_t=x]\,,&\; w(x,t)=\mathbb{E}[\Gamma_{0}|V_{t^-}=x]\,,&\; w_0(x)=\mathbb{E}[\Gamma_{0}|V_0=x]\,,\\
	\Gamma_t=\frac{\partiel V_t}{\partiel t}\,,&\; a(x,t)=\var[\Gamma_{t}|V_t=x]\,.&
\end{eqnarray*}
We recall the conditional variance
$$\var[\Gamma_{t}|V_t=x]=\mathbb{E}[\Gamma_{t}^2|V_t=x]-\left(\mathbb{E}[\Gamma_{t}|V_t=x]\right)^2.$$

\begin{Lemma}\label{th accel_ec_limite} 
For all $t$, we have
\begin{align*}
	\mu^+(\cdot,t)&=\mu(\cdot,t)\;,\quad \mu^-(\cdot,t)=\mu(\cdot,t^-)=\law(V_{t^-})\,,\\
	w(V_{t},t)&=\E[\Gamma_0|V_t]=\E[\Gamma_t|V_t]\,,\\
	w^-(V_{t^-},t)&=\E[\Gamma_0|V_{t^-}]=\E[\Gamma_{t^-}|V_{t^-}]\;,\\
	a(V_t,t)&=\var[\Gamma_{t}|V_t]=\mathbb{E}[\Gamma_{t}^2|V_t]-\left(\mathbb{E}[\Gamma_{t}|V_t]\right)^2\,.
\end{align*}
If $\mu(\cdot,t)=\mu^-(\cdot,t)$, then $w(\cdot,t)=w^-(\cdot,t)$. Moreover, 
we have, in the sense of Radon-Nykodim product of a function and a measure :
\begin{align*}
	(w\mu)^+=w\mu\,,\quad (w\mu)^-=w^-\mu^-\,.
\end{align*}
\end{Lemma}

While all of these claims are readily verifiable, we nonetheless detail the proof of $\E[\Gamma_0|V_{t^-}]=\E[\Gamma_{t^-}|V_{t^-}]$. For any two successive shock times $s_1<s_2=t$, the process $\Gamma$ is constant on $[s_1,t[$; so $$\Gamma_{t^-}=\Gamma_{s_1}=\E[\Gamma_0|X_{s_1}]\,.$$
Moreover,  the r.v. $V_{t^-}$ is  $\sigma(X_{s_1})$-measurable.  So for any Borel bounded function $f$,
$$\E[f(V_{t^-})\Gamma_0]=\E[f(V_{t^-})\E[\Gamma_0|X_{s_1}]]=\E[f(V_{t^-})\Gamma_{t^-}]\,.$$

\begin{proposition}\label{th gas disc}
We have 
	\begin{align}\label{gsp_vitesse}
		\begin{cases}
			\partial_t(\mu) + \partial_x(w\mu) = \Delta(\mu)\\
			\partial_t(w\mu) + \partial_x((w^2+a)\mu) = \Delta(w\mu)\\
			\mu|_{t=0}=\mu_0\,,\quad w\mu|_{t=0}=w_0\mu_0\,,\quad  a\mu|_{t=0}=0\\
			\Delta(\mu)|_{t=0}=0\,,\quad  \Delta(w\mu)|_{t=0}=0\,.
		\end{cases}
	\end{align}
\end{proposition}
Note, particularly in this case, that $\Delta(\mu)$ and $\Delta(w\mu)$ are finite sums of measures, since they are in fact indexed on shock times (the potential discontinuity times of the velocity function).
It is worth noting that it is the shift from continuous measures (of (\ref{Eq2})) to discontinuous $\mu$ and $w\mu$ which gives rise to the \textit{jump} measures $\Delta(\mu)$ and $\Delta(w\mu)$. If $\mu$ and $w\mu$ where continuous, then the right-hand side of the system would be null. So the rise of their \textit{jumps} in the equations is a justified trade-off for transitioning to the setting of discontinuous measures. That's why we do not consider them as forces. Another significant aspect is the appearance of the function $a$, which accounts for the variance of acceleration at a fixed velocity. 
For now, let us look forward to the \textit{jump} measures to gain insight.
\begin{Remark}
	Let us define the stochastic measure on $\R\times\R_+^*$ :
		$$m=\sum_{t\in\R_+^*}(\delta_{V_{t^+}}-\delta_{V_{t^-}})\delta_t\,,$$
	where $\delta_x$ stands for the Dirac measure at $x$. For any bounded Borel function $f$,
	\begin{align*}
		\int_{\R_+^*}\int_{\R} f\partiel \Delta(\mu)=\E\left[\int f\partiel m\right]=\sum_{t\in\R_+^*}\E[f(V_{t^+})-f(V_{t^-})]\,,\\
		\int_{\R_+^*}\int_{\R} f\partiel \Delta(w\mu)=\E\left[\Gamma_0\int f\partiel m\right]
		=\sum_{t\in\R_+^*}\E\left[\Gamma_0(f(V_{t^+})-f(V_{t^-}))\right]\,.
	\end{align*}
\end{Remark}
Note that $m$ is well defined since the sum is in fact indexed on shock times (the potential discontinuity times of the velocity function).
The signed measure $m$ tracks changes of the velocity function over time, especially at discontinuities. Indeed, for any real subset $A$,
$$(\delta_{V_{t^+}}-\delta_{V_{t^-}})(A)=\left\{
\begin{array}{rl}
	+1 &\mbox{if the process $V$ enters the set $A$ at time $t$;}\\
-1 &\mbox{if it exits  the set $A$ at time $t$;}\\
0  &\mbox{if it remains inside or outside $A$.}
\end{array}\right.
$$
Then, we can interpret $m$ as a stochastic threshold-crossing measure of the velocities. So the measures $\Delta(\mu)$ and $\Delta(w\mu)$ can be interpreted as average threshold-crossing metrics of the velocities.  \\

\noindent\textbf{Proof of Proposition.}\\
For the first equation, consider any $f\in C_{0}^{1}(\mathbb{R})$ and $0<t_1<t_2$. We have
	\begin{align*}
	\Delta:&=\int f(x)\mu(dx,t_2) - \int f(x)\mu(dx,t_1)= \mathbb{E}[f(V_{t_2})-f(V_{t_1})]\\
	&= \mathbb{E}\left[\int_{t_1}^{t_2}f'(V_t)\Gamma_tdt+\sum_{t_1< t\leq t_2}(f(V_{t^+})-f(V_{t^-}))\right]\,.
	\end{align*}
Since $V_t$ is $\sigma(X_t)$-measurable and $\Gamma_t=\mathbb{E}[\Gamma_{0}|X_t]$, we get
	\begin{align*}
		\Delta&= \int_{t_1}^{t_2}\mathbb{E}\left[f'(V_t)\Gamma_0\right]dt
		+\sum_{t_1< t\leq t_2}\E[f(V_{t^+})-f(V_{t^-})]\\
	&= \int_{t_1}^{t_2}\mathbb{E}\left[f'(V_t)w(V_t,t)\right]dt
	+\sum_{t_1< t\leq t_2}\int f(x)\partiel[\mu(x,t^+)-\mu(x,t^-)]\\
	&= \int_{t_1}^{t_2}\int f'(x)w(x,t)d\mu(x,t)dt+\sum_{t_1< t\leq t_2}\int f(x)\partiel[\mu(x,t^+)-\mu(x,t^-)]\;.
\end{align*}
For the second equation, we have

\begin{align*}
	\Delta:&=\int f(x)w(x,t_2)d\mu(x,t_2)-\int f(x)w(x,t_1)d\mu(x,t_1)\\
	&=\mathbb{E}[f(V_{t_2})w(V_{t_2},t_2)]-\mathbb{E}[f(V_{t_1})w(V_{t_1},t_1)]
	=\mathbb{E}[\Gamma_0(f(V_{t_2})-f(V_{t_1}))]\\
	&= \mathbb{E}\left[\Gamma_0\int_{t_1}^{t_2}f'(V_t)\Gamma_tdt+\sum_{t_1< t\leq t_2}\Gamma_0(f(V_{t^+})-f(V_{t^-}))\right]\\
	&=\int_{t_1}^{t_2}\mathbb{E}\left[f'(V_t)\Gamma_0\Gamma_t\right]dt
	+\sum_{t_1< t\leq t_2}\Big(\E[f(V_{t^+})\Gamma_0]-\E[f(V_{t^-})\Gamma_0]\Big)\\
	&=\int_{t_1}^{t_2}\mathbb{E}[f'(V_t)\mathbb{E}[\Gamma_0|X_t]\Gamma_t]dt
	+\sum_{t_1< t\leq t_2}\Big(\E[f(V_{t^+})\E[\Gamma_0|V_{t^+}]]-\E[f(V_{t^-})\E[\Gamma_0|V_{t^-}]]\Big)\;.
	\end{align*}
	Applying Lemma \ref{th accel_ec_limite}, we get
	
	\begin{align*}
	\Delta
	&=\int_{t_1}^{t_2}\mathbb{E}[f'(V_t)\Gamma_t^2]dt
	+\sum_{t_1< t\leq t_2}\Big(\E[f(V_{t^+})w^+(V_{t^+},t)]-\E[f(V_{t^-})w^-(V_{t^-},t)]\Big)\\
	&=\int_{t_1}^{t_2}\mathbb{E}[f'(V_t)w(V_t,t)^2]dt
	+\int_{t_1}^{t_2}\mathbb{E}[f'(V_t)\var(\Gamma_t|V_t)]dt\\
	&+\sum_{t_1< t\leq t_2}\int f(x)\partiel[w^+(x,t)\mu(x,t^+)-w^-(x,t)\mu(x,t^-)]\\
	&=\int_{t_1}^{t_2}\int f'(x)w(x,t)^2d\mu(x,t)dt
	+\int_{t_1}^{t_2}\int f'(x)a(x,t)d\mu(x,t)dt\\
	&+\sum_{t_1< t\leq t_2}\int f(x)\partiel[w^+(x,t)\mu(x,t^+)-w^-(x,t)\mu(x,t^-)]\;.
\end{align*}

The weak convergence of $(\mu,w\mu)$, as $t$ tends to zero, is  shown in the same way as for $(\rho,u\rho)$ in (\ref{Eq2}). 

For weak convergence of $a\mu$ towards zero, consider the first shock time $T$ and define 
$\delta=\inf\{t>0 : \exists (i,j)\; s.t.\; v_i\neq v_j,\;v_i+t\theta_i=v_j+t\theta_j\}\wedge T$. We have 
$$a(v_i(t),t)=0\,,\quad \forall\, i,\forall\,t\in]0,\delta[\,.$$
So $(a\mu)(\cdot,t)=0$ for all $t\in]0,\delta[$.\\
To end the proof, remark that $V_{t^+}=V_{t^-}$ for all $t\in]0,T[$. So 
$\mu(\cdot,t^+)=\mu(\cdot,t^-)$ and $w^+(\cdot,t)=w^-(\cdot,t)$ for all $t\in]0,T[$.
\hfill$\square$\\

Let us look forward to the congestion-like term.

\begin{Remark}
	\begin{enumerate}
		\item $a(v,t)=0$ means that all the particles (clusters) of velocity $v$ at time $t$ have the same acceleration. This implies that in the next collision event involving one of these clusters, the others can not participate in the same interaction.
		\item There exists $\delta>0$ such that $a(\cdot,t)\equiv0$ for all $t\in]0,\delta[$. This means that there is a delay before the onset of congestion-like effects. 
		\item Suppose that $a_0:=a(\cdot,0)\not\equiv0$. We have  
			$$\lim\limits_{t\to0}(a\mu)(\cdot,t)=0\neq a_0\mu_0\,.$$ 			
	\end{enumerate}
\end{Remark}

\begin{Corollary}
	\begin{enumerate}
		\item $(\mu,w\mu,0)$ is a weak solution of the gas system (\ref{Eq2}) if and only if :
		\begin{enumerate}
			\item $\mu$ is continuous; 
			\item $a=0$, $\partiel\mu\partiel t$-p.p. This is equivalent to
				$$\Gamma_t=\mathbb{E}[\Gamma_t|V_t]\;\mbox{a.s. and for almost all } t,$$
			which means that acceleration is a function of velocity.
		\end{enumerate}
		\item Assume that particles with the same initial velocities also have the same initial accelerations. If $T$ describes the first shock time, then $(\mu,w\mu,0)$ is a weak solution of the gas system (\ref{Eq2}) on $\R\times]0,T[$. 
	\end{enumerate}
\end{Corollary}
In future research, it would be worthwhile to investigate the initial conditions that guarantee the continuity of $\mu$ and the vanishing of the conditional variance $a$.

\section*{Conclusion and perspectives}
We have introduced a generalized variational principle for a discrete sticky particle model subjected to external forces. We have also introduced an innovative solution to the unforced pressureless gas system with a congestion-like term, relying on the particle acceleration function and the mass distribution over particle velocities. Unlike traditional approaches based on particle positions, this velocity-space perspective enables the exploration of regimes where particles are characterized by their kinematic state rather than their position. This approach opens up interesting possibilities for the analysis of kinematic data and modeling in statistical physics, while providing an extension to discontinuous and non-discrete measures.

{}
\end{document}